\documentstyle[12pt,amssymb]{article}
\input xy
\xyoption{all} \hfuzz=1.5pt \tolerance=500

\newcommand{\mHtp}{\mathbin{\stackrel{\cdot}{\otimes}}}
\newcommand{\mH}{\mHtp}
\newcommand{\mh}{\mH}
\newcommand{\ot}{\otimes}

\newcommand{\tr}{\triangleright}
\newcommand{\tl}{\triangleleft}

\newcommand{\va}{\varphi}

\newcommand{\id}{{\bf 1}}
\newcommand{\co}{{\B C}}

\newcommand{\SS}{{\cal S}}
\newcommand{\X}{{\cal X}}
\newcommand{\bd}{\begin{document}}
\newcommand{\ed}{\end{document}}
\newcommand{\Htp}{\mathbin{\stackrel{\cdot}{\bigotimes}}}
\newcommand{\Hts}{\mathbin{\stackrel{\cdot}{\bigoplus}}}
\newcommand{\nl}{\nu\in\Lambda}
\newcommand{\xl}{X_\nu;\nl}
\newcommand{\x}{X_\nu}
\newcommand{\la}{\langle}
\newcommand{\ra}{\rangle}
\newcommand{\ch}{\cal H}
\newcommand{\su}{\subseteq}
\newcommand{\B}{\Bbb}
\newcommand{\e}{\varepsilon}
\newcommand{\vk}{\varkappa}
\newcommand{\vt}{\vartheta}
\newcommand{\al}{\alpha}
\newcommand{\de}{\delta}
\newcommand{\lo}{\Longleftrightarrow}
\newcommand{\k}{{\cal F}}
\newcommand{\cK}{{\cal F}}
\newcommand{\bb}{{\cal B}}
\newcommand{\dd}{{\cal D}}
\newcommand{\r}{{\cal R}}
\newcommand{\cR}{{\cal R}}
\newcommand{\f}{{\cal K}}
\newcommand{\cb}{{\cal CB}}
\newcommand{\wrr}{\widetilde{\cal R}}
\newcommand{\br}{^\bullet{\cal R}}
\newcommand{\en}{E_\nu}
\newcommand{\el}{{\cal L}}
\newcommand{\lr}{\Longrightarrow}
\newcommand{\Long}{\Longleftarrow}
\newcommand{\q}{\quad}
\newcommand{\qq}{\qquad}
\newcommand{\cd}{\cdot}
\newcommand{\fur}{f: E\to{\B R}}
\newcommand{\furo}{f_0: E_0\to{\B R}}
\newcommand{\fuc}{f: E\to{\B C}}
\newcommand{\ii}{\infty}
\newcommand{\di}{\diamondsuit}
\newcommand{\bgd}{{\bigtriangledown}}
\newcommand{\bu}{{\bigtriangleup}}
\newcommand{\bc}{{completely bounded}}
\newcommand{\cc}{{completely contractive}}
\newcommand{\qs}{{quantum space}}
\newcommand{\isc}{{isometric}}
\newcommand{\ism}{{isomorphism}}
\newcommand{\qss}{{quantum spaces}}
\newcommand{\bco}{{completely bounded operator}}
\newcommand{\bcos}{{completely bounded operators}}
\newcommand{\res}{{respectively}}
\newcommand{\tp}{{tensor product}}
\newcommand{\eq}{{equivalent}}
\newcommand{\qtp}{{quantum tensor product}}
\newcommand{\mma}{\mathrel{\mathop{\otimes}\limits_{A}}}
\newcommand{\mmb}{\mathrel{\mathop{\otimes}\limits_{\bb}}}
\newcommand{\mmp}{\mathrel{\mathop{\otimes}\limits_{p}}}
\newcommand{\mmh}{\mathrel{\mathop{\otimes}\limits_{h}}}
\newcommand{\mmf}{\mathrel{\mathop{\otimes}\limits_{4}}}
\newcommand{\mmi}{\mathrel{\mathop{\otimes}\limits_{i}}}
\newcommand{\mmm}{\mathrel{\mathop{\otimes}\limits_{\bb-\bb}}}
\newcommand{\mmo}{\mathrel{\mathop{\cdot}\limits_{1}}}
\newcommand{\mmA}{\mathrel{\mathop{\otimes}\limits_{A-A}}}
\newcommand{\mmd}{\mathrel{\mathop{\cdot}\limits_{2}}}
\newcommand{\msp}{\mathrel{\mathop{\otimes}\limits_{sp}}}
\newcommand{\mms}{\stackrel{h}{\otimes}}
\newcommand{\mmt}{\stackrel{4}{\otimes}}
\newcommand{\mmx}{\stackrel{i}{\otimes}}
\newcommand{\mmy}{\stackrel{p}{\otimes}}
\newcommand{\mmz}{\stackrel{sp}{\otimes}}
\newcommand{\gd}{\ddagger}
\newcommand{\od}{\odot}
\newcommand{\mt}{\mapsto}
\newcommand{\mmc}{\mathrel{\mathop{\otimes}\limits_{\sim}}}
\newcommand{\mme}{\stackrel{\sim}{\otimes}}

\begin{document}

   \centerline{{ \bf Extreme flatness of normed modules and}}

   \centerline{{\bf Arveson-Wittstock type theorems}\footnote{This research
was supported by the Russian Foundation for Basic Research (grant
No. 05-01-00982).}}

\bigskip

   \centerline{\it Dedicated to the memory of Graham R. Allan}

   \vspace{1cm}

\centerline{A.~Ya.~Helemskii}

\centerline{Faculty of Mechanics and Mathematics}

\centerline{Moscow State University}

\centerline{Moscow 119992 Russia}

\bigskip

\noindent We show in this paper that a certain class of normed
modules over the algebra of all bounded operators on a Hilbert space
possesses a homological property which is a kind of a
functional-analytic version of the standard algebraic property of
flatness. We mean the preservation, under projective tensor
multiplication of modules, of the property of a given morphism to be
isometric. As an application, we obtain several extension theorems
for different types of modules, called Arveson--Wittstock type
theorems. These, in their turn, have, as a straight corollary, the
`genuine' Arveson-Wittstock Theorem
in its non-matricial presentation. We recall that the latter theorem
plays the role of a  `quantum'  version of the classical
Hahn--Banach theorem on the extension of bounded linear functionals.
It was originally proved in~\cite{wit}, and a crucial preparatory
step was done in~\cite{ar}. As to the monographical presentation,
see the textbooks~\cite{{ER},{pu}}.

\bigskip
{\bf Acknowledgment.} We would like to express our profound
gratitude to the referee for his tremendous job. Apart from all
things that are usually required from a competent report, he has
shown to the author, a novice in the field, how to put his results
in the right perspective, suggesting, in particular, the lines of
their future continuation. More of this, with the help of a rather
sophisticated argument he has obtained statements, augmenting and
strengthening some central results of our paper. Finally, he drew
our attention to a number of valuable articles, related to the
subject. A part of his contribution is reflected in Remarks 4 and 6
below.

\bigskip

\centerline{\bf $0$. Preliminaries}

\bigskip

\noindent Throughout the paper we shall denote by ${\cal B}(E,F)$
the space of all bounded operators  acting between normed spaces $E$
and $F$,
 always
equipped with the operator norm. We shall denote by ${\cal F}(E,F)$
the
  subspace of this space consisting of the finite-rank operators.
As usual, we set  ${\cal B}(E):={\cal B}(E,E)$ and ${\cal
F}(E):={\cal F}(E,E)$.

The identity operator on $E$ is denoted by $\id_E$.

The inner product in Hilbert spaces is denoted by
$\la\,\cd\,,\,\cd\,\ra$. The complex-conjugate Hilbert space of a
Hilbert space $H$ is denoted by $H^c$.

In our future arguments we shall come across some tuples of
isometric operators between Hilbert spaces, say $H$ and $K$. Let
$S_k;k=1,...,n$ be such a tuple, and suppose that the final
projections $P_k:=S_kS_k^*$ of these operators are pairwise
orthogonal. We recall that in this situation we have the following
equalities:
$$
S_k=P_kS_k,\; S_k^*=S_k^*P_k\; {\rm and},{\rm{\; as\; a\;
corollary}},\; S_k^*S_l=0 \; {\rm for}\; k\ne l.
 \eqno(1)
$$

Another class of operators we shall need is that of rank-one
operators. For the same $H$ and $K$ as above,
 and for $\xi\in K$ and $\eta\in H$, we denote by $\xi\bigcirc\eta$ the rank-one operator taking $\zeta\in H$ to
$\la\zeta,\eta\ra\xi\in K$. We recall the equalities
$$
(\xi\bigcirc\eta)(\xi'\bigcirc\eta')=\la\xi',\eta\ra \xi\bigcirc
\eta',\; a(\xi\bigcirc\eta)= (a\xi)\bigcirc\eta\; {\rm and}\;
(\xi\bigcirc\eta)a=\xi\bigcirc(a^*\eta),\eqno (2)
$$
that are valid whenever their ingredients make sense.

As usual, the symbol `$\ot$' denotes the algebraic tensor product of
linear spaces and operators. Further, we use the symbol
 ` $\mHtp$ '  for the Hilbert tensor product of Hilbert spaces as
well as for the Hilbert tensor product of operators acting between
these spaces. Finally, the symbol  ` $\mmp$ '  denotes the {\it
non-completed} projective tensor product of normed spaces.

Further, we choose a separable infinite-dimensional Hilbert space,
denote it by $L$, and fix it throughout the whole paper. Sometimes
in what follows this Hilbert space will be refered as the
`canonical' one.\footnote{Our experience shows that, as a whole, it
is more convenient to make an  `abstract' choice, and not be tied
to, say, $l_2$ or $L^2(\,\cd\,)$.} For brevity, we denote the
operator algebras ${\cal B}(L)$ and $\k(L)$ by ${\cal B}$ and $\k$,
respectively.

Throughout the paper, {\it the terms left module, right module and
 bimodule (=two-sided module) always mean a unital module of the relevant type
over the operator algebra $\bb$}; we shall never consider other
basic algebras.
 The respective outer (=module) multiplications will be denoted by a dot:  `$\,\cd\,$'.
The words  {\it (bi)module morphism}  always mean a morphism of the
$\bb$-(bi)modules in question.

Let $X$ be a left (respectively, right) module. A {\it left}
(respectively, {\it right}) {\it support} of the element $u\in X$
is, by definition, each  projection $P\in{\cal B}$ such that $P\cdot
u=u$ (respectively, $u\cdot P=u$). If we have a bimodule, and $P$ is
both a left and  right support of the  element $u$, then we say that
$P$ is (just) a {\it support} of $u$.

Let $X$  be a left module and simultaneously a normed space. We
recall that $X$ is called a {\it contractive} left module  if we
have $\|a\cd u\|\le\|a\|\|u\|$ for all $a\in\bb$ and $u\in X$.
Similarly, the condition $\|u\cdot a\|\le\|a\|\|u\|$ leads to the
notion of a contractive right module, and the two  mentioned
conditions together lead to the notion of a contractive bimodule.

If $X$ and $Y$ are two contractive left modules, we  denote the
space of all bounded (as operators) morphisms between $X$ and $Y$ by
$_\bb{\bf h}(X,Y)$. The relevant spaces for the cases of right
modules and bimodules will be denoted by ${\bf h}_\bb(X,Y)$ and
$_\bb{\bf h}_\bb(X,Y)$, respectively. We equip these spaces with the
operator norm, that is, we consider them as normed subspaces of
${\bb}(X,Y)$.

Let $X$ be a contractive left module. Then the  complex conjugate
normed space $X^{c}$ becomes a contractive right module with the
outer multiplication $x\cd a$, defined as the former $a^*\cd x$.
Similarly, a right outer multiplication on $X$ gives rise to a left
one on $X^{c}$, defined by $a\cd x:=x\cd a^*$. We call $X^{c}$,
equipped with the relevant structure of the contractive right or
left module, {\it the complex conjugate module of $X$}. Obviously,
every bounded morphism $\va:X\to Y$ of contractive left
(respectively, right) modules, being considered as a map from
$X^{c}$ into $Y^{c}$, becomes a morphism of right (respectively,
left) modules with the same norm.

We recall several standard constructions. Let $X$ be a left
contractive module. Then its dual space $X^*$ is a right contractive
module with the outer multiplication defined by $$[f\cd
a](x):=f(a\cd x); a\in A,x\in X,f\in X^*\,.
$$ Similarly, the dual to a right contractive module becomes a left contractive module
with the help of the equality $[a\cd f](x):=f(x\cd a)$, and the dual
to a contractive bimodule becomes itself a contractive bimodule with
the help of both of these equalities. If $X$ and $Y$ are two left
contractive modules, then the normed space $\bb(X,Y)$ is a
contractive bimodule with outer multiplications defined by
$$[a\cd\va](x):=a\cd(\va(x))\quad {\rm  and} \quad  [\va\cd
a](x):=\va(a\cd x)\,;$$
here $\va\in\bb(X,Y)$, etc. Finally, if $X$ is a left and $Y$ is a
right contractive module, then the normed space $X\mmp Y$ is a
contractive bimodule with the outer multiplications uniquely defined
by
$$
a\cd(x\ot y):=(a\cd x)\ot y \quad {\rm  and}\quad   (x\ot y)\cd
a:=x\ot(y\cd a)\,.
$$

We shall  also need the notion of the module and bimodule {\tp}s, in
their projective non-completed version. Suppose that either $X$ is a
right  and $Y$ is a left contractive modules, or both of $X$ and $Y$
are contractive bimodules. The normed spaces $X\mmb Y$ in the first
case, and $X\mmm Y$ in the second one, called respectively {\it the
module and bimodule {\tp} of $X$ and $Y$}, are defined in terms of
the universal property with respect to the class of balanced,
bounded, bilinear operators from $X\times Y$ into normed spaces;
cf., e.g.,~\cite{he2}. Namely, in the   one-sided   case a bounded
bilinear operator $\r:X\times Y\to E$, where $E$ is a normed space,
is called {\it balanced} if $ \r(x\cd a,y)=\r(x,a\cd y)$ for all
$a\in \bb$, $x\in X$, and $y\in Y$. In the two-sided case such a
bilinear operator is called {\it balanced}, if, in addition to
  the indicated equalities, we also have $\r(a\cd x,y)=\r(x,y\cd a)$.

As to  explicit constructions, the spaces $X\mmb Y$ and $X\mmm Y$
can be realized as the normed quotient spaces of $X\mmp Y$, the
projective tensor product of the underlying normed spaces of our
(bi)modules. Namely, $$X\mmb Y=X\mmp Y/N_1\,,
$$
 where $N_1$ is the closure of $${\rm span}\{ x\cd a\ot y-x\ot a\cd y;
 a\in\bb, x\in X, y\in Y\}\,,
$$
 whereas
$$
X\mmm Y=X\mmp Y/N_2\,,
$$
 where $N_2$ is the closure of $${\rm span}\{a\cd x\ot y-x\ot y\cd a,
x\cd a\ot y-x\ot a\cd y; a\in\bb, x\in X, y\in Y\}\,.$$
Consequently, for the elementary tensors in
 $X\mmb Y$ and $X\mmm Y$, that is, cosets $x\mmb y:=x\ot y+N_1$ and $x\mmm y:=x\ot y+N_2$,
we have the identities
$$
x\cd a\mmb y=x\mmb a\cd y,\q a\cd x\mmm y=x\mmm y\cd a\q, {\rm
and}\q x\cd a\mmm y=x\mmm a\cd y.\eqno(3)
$$
Finally, the norm of an element $u$ in $X\mmb Y$ or in $X\mmm Y$ is
equal to
$$
\inf\left\{\sum_{k=1}^n\|x_k\|\|y_k\|\right\},\eqno(4)
$$
\noindent where the infimum is taken over all possible
representations of $u$ in the form $\sum_{k=1}^nx_k\mmb y_k$ or,
according to the case, $\sum_{k=1}^nx_k\mmm y_k$.

Note the following attractive property of module {\tp}s over $\bb$.

\medskip
{\bf Proposition 1.} {\it Let $X$ be a right  and   $Y$ a left
module. Then every $u\in X\mmb Y$ can be represented as a single
elementary tensor.
Moreover, if $$u=\\
\sum_{k=1}^nx_k\mmb y_k;x_k\in X,y_k\in Y\,,$$ and $S_k;k=1,...,n$
is an arbitrary family of {\isc} operators on $L$ with pairwise
orthogonal final projections
 $P_k:=S_kS_k^*$, then such a representation can be taken as $u=x\mmb y$, where
$x:=\sum_{k=1}^nx_k\cd S^*_k$ and
 $y:=\sum_{k=1}^nS_k\cd y_k$.}

\smallskip
$\tl$ By (1) and (3), we have
$$
x\mmb y=\sum_{k,l=1}^nx_k\cd S^*_k\mmb S_l\cd
y_l=\sum_{k,l=1}^nx_k\cd S^*_kS_l\mmb \cd y_l= \sum_{k=1}^nx_k\mmb
y_k=u. \q\tr
$$

\medskip
Finally, we recall that the construction of the module tensor
product has functorial properties. Namely, if $\al:X_1\to X_2$ and
$\beta:Y_1\to Y_2$ are bounded morphisms of contractive right and
left modules, respectively, then there exists a bounded operator
$\al\mmb\beta:X_1\mmb Y_1\to X_2\mmb Y_2$, uniquely defined by the
rule $x\mmb y\mt\al(x)\mmb\beta(y)$. Moreover, we have
$\|\al\mmb\beta\|\le\|\al\|\|\beta\|$. If we deal with contractive
bimodules and their bimodule morphisms, then there exists a bounded
operator $\al\mmm\beta:X_1\mmm Y_1\to X_2\mmm Y_2$, uniquely defined
by the rule $x\mmm y\mt\al(x)\mmm\beta(y)$, and we have
$\|\al\mmm\beta\|\le\|\al\|\|\beta\|$.

\bigskip

\centerline{\bf $1$. Ruan bimodules and semi-Ruan one-sided modules}

\bigskip

\medskip
{\bf Definition 1} (cf. the definition of an  `abstract operator
space'  in~\cite[p.\ 20]{ER} or~\cite[p.\ 180-181]{pu}). A
contractive bimodule $Y$ is called a {\it Ruan bimodule} if it
satisfies the following condition (the {\it Ruan axiom}):

\q $(R)$: for all $u,v\in X$ with orthogonal supports, we have
$$
\|u+v\|=\max\{\|u\|,\|v\|\}.
$$


\medskip
{\bf Example 1.} Consider the Banach space $\bb(L\mh H,L\mh K)$,
where $H$ and $K$ are arbitrary Hilbert spaces (whatever Hilbert
dimension, finite or infinite, they would have). It is easy to check
that this space is a Ruan bimodule with respect to the outer
multiplications, defined by $a\cd\tilde b:=(a\mh\id_K)\tilde b$ and
$\tilde b\cd a:=\tilde b(a\mh\id_H); a\in\bb,\tilde b\in\bb(L\mh
H,L\mh K)$.

\medskip
For Ruan bimodules, the axiom $(R)$ can be strengthened.

\medskip
{\bf Proposition 2}. {\it Let $u_1,\dots,u_n$ be elements of a Ruan
bimodule $X$ with pairwise orthogonal left supports, say $P_k$, and
pairwise orthogonal right supports, say $Q_k;k=1,\dots,n$. Then}
$$
\|u_1+\cdots+u_n\|=\max\{\|u_1\|,\dots,\|u_n\|\}.
$$

\smallskip
$\tl$ For  brevity, set $u:=\sum_{k=1}^nu_k$, and take an arbitrary
tuple $S_k;k=1,...,n$ of isometric operators on $L$ with pairwise
orthogonal final projections.

At first we compare the norms of the elements $u$ and  $v:=
\sum_{k=1}^nS_k\cd u_k\cd S_k^*$. Set also $a:=\sum_{k=1}^nP_kS_k^*$
and $b:=\sum_{k=1}^nS_kQ_k$. Then the equalities (1)
imply that
\begin{eqnarray*}
a\cd v\cd b&=&\sum_{k,l,m=1}^nP_kS_k^*\cd(S_l\cd u_l\cd S_l^*)\cd S_mQ_m\\
&=& \sum_{k,l,m=1}^nP_k\cd(S_k^*S_l\cd u_l\cd S_l^*S_m)\cd
Q_m=\sum_{k=1}^nP_k\cd u_k\cd Q_m=u\,.
\end{eqnarray*}
Therefore $\|u\|\le\|a\|\|v\|\|b\|$. Further, the $C^*$-identity
gives $\|a\|=\|aa^*\|^ {1/2}$ and $\|b\|=\|b^*b\|^{1/2}$. Again
using   (1), we have
$$aa^*=\sum_{k,l=1}^nP_kS_k^*S_lP_l^*=\sum_{k=1}^nP_k\,,$$  and similarly $b^*b=\sum_{k=1}^nQ_k$. Thus
both of these operators are projections, and hence their norm is 1.
Consequently, $\|u\|\le\|v\|$.

Now observe that, by the same equalities (1), the final projection
of $S_k$ is a support of the element $S_k\cd u_k\cd S_k^*\in X;
k=1,\dots,n$. Since these projections are pairwise orthogonal, it
follows from $(R)$ that
$$
\|v\|= \max\{\|S_1\cd u_1\cd S_1^*\|,\dots,\|S_n\cd u_n\cd
S_n^*\|\}\,.
$$
 Since $X$ is contractive, and
$S_k^*S_k=\id_L$, we have, for every $k$, $\|S_k\cd u_k\cd
S_k^*\|=\|u_k\|$. Therefore $\|v\|=\max\{\|u_1\|,\dots,\|u_n\|\}.$

Thus $\|u\|\le\max\{\|u_1\|,\dots,\|u_n\|\}$. But we obviously have
$u_k=P_k\cd u\cd Q_k$, and our bimodule is contractive. From this,
we have the reverse inequality. $\tr$

\medskip
We turn from bimodules to one-sided modules. As  experience shows,
the obvious version of the condition (R) for these modules is not
very workable. The following, more `tolerant' definition happens to
be more useful.

\medskip
{\bf Definition 2.} A contractive left module $X$ is  a {\it left
semi-Ruan module}\footnote{B.Magajna in~\cite[Corollary 2.2]{mag},
pursuing different aims, considers a certain class of left modules
over arbitrary $C^*$-algebras. It is not hard to see that in the
case when the algebra in question is $\bb$, this class coincides
with the class of Banach semi-Ruan modules. We are indebted to
D.Blecher, who drew our attention to the paper of Magajna.}, if it
satisfies the following condition:

$(lsR)$: if $u,v\in X$ have orthogonal left supports, then
$$
\|u+v\|\le(\|u\|^2+\|v\|^2)^{1/2}.
$$

Similarly, with the obvious modifications, we introduce the notion
of a {\it right semi-Ruan module}. The respective condition will be
denoted by $(rsR)$.

Needless to say, we have a similar estimate for several summands.
Namely, {\it if elements $u_1,\cdots,u_n$ of a one-sided semi-Ruan
module have respective one-sided pairwise orthogonal supports, then
$\|u_1+\cdots + u_n\|\le(\|u_1\|^2+\cdots+\|u_n\|^2)^{{1/2}}$.}

Clearly, every sub-bimodule of a Ruan module is itself a Ruan
module, and   similar hereditary property holds for one-sided
semi-Ruan modules. Note also the following obvious observation.

\medskip
{\bf Proposition 3.} {\it The complex conjugate module} (cf. the
previous section) {\it of a semi-Ruan module is itself a semi-Ruan
module.} $\tl \tr$

\medskip
Here is our most important pair of examples.

\medskip
{\bf Example 2.} For an arbitrary Hilbert space $H$, the Hilbert
space $L\mh H$ is obviously a left semi-Ruan module with respect to
the outer multiplication
$$
a\cd\zeta:=(a\mh\id_H)\zeta; a\in\bb, \zeta\in L\mh H\,.
$$
Its complex conjugate right semi-Ruan module  is, of course, the
Hilbert space $L^{c}\mh H^{c}$ with the outer multiplication
$\zeta\cd a:=(a^*\mh\id_H)\zeta$. Note that this latter module is,
by virtue of the Riesz representation theorem, nothing else than the
dual to the left module $L\mh H$.

\medskip
{\bf Proposition 4} (cf.~\cite[Proposition 2]{he4}). {\it Every Ruan
bimodule, considered as a left or right module, is a respective
one-sided semi-Ruan module.}

\smallskip
$\tl$ Let $X$ be our bimodule, and let $u_1,u_2\in X$ have, to be
definite, pairwise orthogonal left supports $P_1$ and $P_2$. Of
course, we may  suppose  that $u_1,u_2\ne0$. Take isometric
$S_1,S_2\in\bb$ with orthogonal final projections, and set
$$
v:=\frac{1}{\|u_1\|}u_1\cd S_1^*+\frac{1}{\|u_2\|}u_2\cd
S_2^*\,,\quad
 P:=P_1+P_2\,,\quad{\rm and}\quad  b:=\|u_1\|S_1+
\|u_2\|S_2\,.
$$
 Then the equalities (1) easily imply that
$$
P\cd v\cd b=u_1+u_2.
$$
Therefore $\|u_1+u_2\|\le\|P\|\|v\|\|b\|$. But, of course, $P$ is  a
projection, and the $C^*$-identity immediately gives
$$\|b\|=(\|u_1\|^2+\|u_2\|^2)^{{1/2}}.$$ Finally, the summands in $v$
obviously have orthogonal left and orthogonal right supports.
Therefore, since $X$ is contractive, Proposition 2 gives $\|v\|=1$.
The rest is clear. $\tr$

\medskip
{\bf Remark 1.} However, a contractive bimodule, which is a left and
a right semi-Ruan module, is not, generally speaking, a Ruan
bimodule. One can take, as a counter-example, $L\mh L^{c}$ or $L\mmp
L^{c}$.

\medskip
Note also, that the $l_2$-sum of a family of one-sided semi-Ruan
modules is also a semi-Ruan module of the same type.

\medskip
{\bf Proposition 5}. {\it Let $X$ be a right semi-Ruan module, $Y$ a
left semi-Ruan module, and $u\in X\mmb Y$. Then
$$
\|u\|= \inf\{\|x\|\|y\|\},
$$
where the infimum is taken over all possible representations of $u$
in the form $u=x\mmb y;x\in X,y\in Y$.}
 (Such representations exist by Proposition 1).

\medskip
$\tl$ Denote the indicated infimum by $\|u\|'$. It follows from (4)
that $\|u\|\le\|u\|'$. Our task is to establish the reverse
inequality.

Take an arbitrary representation of $u$ in the form
$\sum_{k=1}^nx_k\mmb y_k$. Obviously, without loss of generality we
may suppose  that $\|x_k\|=\|y_k\|; k=1,\dots,n$. Let
$S_k,P_k;k=1,...,n,\; x$ and $y$ be as in Proposition 1.
The formulae (1) imply that $P_k$ is the right support of $x_k\cd
S^*_k$ and the left support of $S_k\cd y_k;k=1,\dots,n$. Therefore
the conditions $(rsR)$ and $(lsR)$ imply for our contractive modules
that
 \begin{eqnarray*}
\|u\|'\le\|x\|\|y\|&\le&
\left(\sum_{k=1}^n\|x_k\cd S^*_k\|^2\right)^{{1/2}}\left(\sum_{k=1}^n\|S_k\cd y_k\|^2\right)^{{1/2}}\\
 &=&
 \left(\sum_{k=1}^n\|x_k\|^2\right)^{{1/2}}\left(\sum_{k=1}^n\|y_k\|^2\right)^{{1/2}}\\
&=& \sum_{k=1}^n\|x_k\|^2=\sum_{k=1}^n\|x_k\|\|y_k\|.
\end{eqnarray*}
Taking all possible representations of $u$ as sums of elementary
tensors and using (4), we obtain
  $\|u\|'\le\|u\|$. $\tr$

\medskip
Let $X$ be a contractive bimodule, and let $Y$ be a contractive left
module. We  consider  the space $X\mmb Y$, where $X$ is considered
as a right contractive module. Recall that in this situation $X\mmb
Y$ has the structure of a contractive left module with the outer
multiplications uniquely defined by $a\cd(x\mmb y):=(a\cd x)\mmb y
$. Similarly, if $X$ is a contractive module  and $Y$ is a
contractive bimodule, then the space $X\mmb Y$, where now $Y$ is
considered as a contractive left module, is a contractive right
module with the outer multiplications uniquely defined by $(x\mmb
y)\cd a:=x\mmb(y\cd a)$.

\medskip
{\bf Proposition 6.} {\it Let  $X$ be  a Ruan bimodule, and let $Y$
be  a left semi-Ruan module. Then $X\mmb Y$ is a left semi-Ruan
module.

Let  $X$ be a right semi-Ruan module, and let $Y$ be a Ruan
bimodule. Then $X\mmb Y$ is a right semi-Ruan module.}

\smallskip
$\tl$ Since the arguments concerning both assertions are strictly
parallel, we restrict ourselves to the first one. Since $X$ is
contractive as a left module, the equality (4) obviously implies
that $X\mmb Y$ is also contractive as a left module. So we
concentrate on the condition $(lsR)$.

Let $u_1,u_2\in X\mmb Y$ have orthogonal left supports, say $Q_1$
and $Q_2$. By virtue of Proposition 1, we may  suppose that
$u_k=x_k\mmb y_k:k=1,2$. Obviously, without loss of generality we
may also suppose  that $\|x_k\|=1$ and $x_k:=Q_k\cd x_k;k=1,2$.

Take, for our $u_k$, the operators $S_k$ and $P_k;k=1,2$ as in the
just-mentioned proposition. Then we have $$u_1+u_2=(x_1\cd
S_1^*+x_2\cd S_2^*)\mmb(S_1\cd y_1+S_2\cd y_2)\,.
$$
 Further, the elements
$x_k\cd S_k^*;k=1,2$ have orthogonal left supports $Q_k$ and
orthogonal right supports $P_k$, respectively. Therefore, since $X$
is a contractive bimodule, Proposition 2 implies that
$$\|x_1\cd S_1^*+x_2\cd S_2^*\|=\max\{\|x_1\cd S_1^*\|,\|x_2\cd S_2^*\|\}=
\max\{\|x_1\|,\|x_2\|\}=1\,.
$$ Consequently, $\|u_1+u_2\|\le\|S_1\cd y_1+S_2\cd y_2\|$. But the elements
$S_k\cd y_k;k=1,2$ have orthogonal left supports $P_k$, and $Y$ is
contractive and satisfies $(lsR)$. Thus
\begin{eqnarray*}
\|u_1+u_2\|&\le&(\|S_1\cd y_1\|^2+\|S_2\cd y_2\|^2)^{{1/2}}\\
&=&(\|y_1\|^2+\|y_2\|^2)^{{1/2}}=
((\|x_1\|\|y_1\|)^2+(\|x_2\|\|y_2\|)^2)^{{1/2}}.
\end{eqnarray*}
It remains to take all possible representations of $u_1$ and $u_2$
as elementary tensors,  and to apply Proposition 5. $\tr$

\bigskip

\centerline{\bf $2$. Extremely flat and extremely injective
(bi)modules}

\bigskip

\noindent We give the following definition in the spirit of the
well-known definitions of  flat and of  strictly flat Banach module
(\cite [Chapter VII, \S1]{he2}, \cite [Chapter VII, \S1.3]{he3}).

\medskip
{\bf Definition 3.} A contractive left module $X$ is   {\it
extremely flat with respect to semi-Ruan modules} or, for short,
{\it ESR-flat}, if, for every isometric morphism
 $\al:Y\to Z$ of right semi-Ruan
modules, the operator $\al\mmb\id_X:Y\mmb X\to Z\mmb X$ (see the end
of Section 0) is also isometric.\medskip

We  define similarly  the  `right-hand'  version of this notion.

Finally, a contractive bimodule $X$ is   {\it extremely flat with
respect to Ruan bimodules} or, for short, {\it ER-flat}, if, for
every isometric morphism $\al:Y\to Z$ of Ruan bimodules, the
operator $\al\mmm\id_X:Y\mmm X\to Z\mmm X$   is also isometric.

\medskip
{\bf Remark 2.} The word  `extremely'  is chosen because isometric
operators or morphisms are exactly the so-called extreme
monomorphisms in some principal categories of spaces or (bi)modules
in functional analysis (cf., e.g.,~\cite{cml},~\cite[Chapter 0,
\S5]{he1}).

\medskip
As simplest examples, the module $\bb$ is ESR-flat as a left and as
a right contractive module, whereas the bimodule $\bb\mmp\bb$ is an
ER-flat contractive bimodule.  Of course, this is because tensoring
by $\bb$ in the one-sided case  and by $\bb\mmp\bb$ in the two-sided
case does not change a given space. In addition, one can easily show
that $\bb\mmp l_1$ and $(\bb\mmp\bb)\mmp l_1$ are ESR-flat as a
one-sided module and ER-flat as a two-sided module, respectively.
Note, that in these examples,   tensoring by the respective
(bi)module preserve the isometry of morphisms of all given
contractive modules, and not only (semi-)Ruan modules. The
properties of the latter modules will be seen to   be indispensable
when, very soon, we   proceed to other examples, more important for
our aims.

We emphasize that the given definition does not require that our
extremely flat (bi)module is itself a (semi-)Ruan (bi)module.
However, in our principal examples that will be the case.

Let us show that several standard constructions preserve the
property of  extreme flatness.

\medskip
{\bf Proposition 7.} {\it If a left or right contractive module is
ESR-flat, then the same is true for its complex conjugate module.}

\smallskip
$\tl$  To be definite, consider a left ESR-module $X$. Our task is
to prove that, for every {\isc} morphism of left semi-Ruan modules
$\al:Y\to Z$, the operator $\id_{X^{c}}\mmb\al:X^{c}\mmb Y\to
X^{c}\mmb Z$ is {\isc}.

Consider $\id_{X^{c}}\mmb\al$ as acting between the respective
complex conjugate normed spaces $(X^{c}\mmb Y)^{c}$ and $(X^{c}\mmb
Z)^{c}$. It is obvious that the first space
 coincides with $Y^c\mmb X$ up to an {\isc} {\ism}, uniquely defined by taking $x\ot y$ to $y\ot x$,
and similarly that the second space
coincides with $Z^c\mmb X$. Moreover, under such an identification
the operator $\id_{X^{c}}\mmb\al$ transforms to
$\al\mmb\id_X:Y^c\mmb X\to Z^c\mmb X$, where
 $\al$, now being considered as a map between $Y^c$ and $Z^c$, is, of course, an {\isc} morphism of the
respective complex conjugate right modules. But the latter are, by
Proposition 3, semi-Ruan modules.

The rest is clear. $\tr$

\medskip
{\bf Proposition 8.} {\it Let $X$ be a left  and $Y$ a right
ESR-flat contractive module. Suppose that at least one of them is a
semi-Ruan module. Then the bimodule $X\mmp Y$} (cf. Section 0) {\it
is ER-flat}.

\smallskip
$\tl$ To be definite, suppose that $Y$ is a semi-Ruan module. Let
$\al:Z_1\to Z_2$ be an {\isc} morphism of Ruan bimodules. Our task
is to show that the operator $\id_{X\mmp Y}\mmm\al: (X\mmp Y)\mmm
Z_1\to(X\mmp Y)\mmm Z_2$ is also {\isc}.

It is known (and easy to verify) that the latter operator is weakly
isometrically equivalent to the operator $(\id_Y\mmb\al)\mmb\id_X:
(Y\mmb Z_1)\mmb X\to(Y\mmb Z_2)\mmb X$. Recall (cf.,
e.g.,~\cite{he1}) that this means that there exists a commutative
diagram
$$
\xymatrix@C+20pt{ (X\mmp Y)\mmm Z_1 \ar[r]^{\id_{X\mmp Y}\mmm\al}
\ar[d]
& (X\mmp Y)\mmm Z_2 \ar[d]\\
(Y\mmb Z_1)\mmb X \ar[r]^{(\id_Y\mmb\al)\mmb\id_X} & {(Y\mmb
Z_2)\mmb X} },
$$
\noindent where the vertical arrows depict isometric isomorphisms of
normed spaces. In the  case that we are considering, these
isomorphisms, a kind of  `complicated associativity', are uniquely
defined by taking an elementary tensor $(x\ot y)\mmm z$ to $(y\mmb
z)\mmb x$; here $x\in X,y\in Y$, and $z$ belongs to $Z_1$ or $Z_2.$

We see that it is sufficient to show that the operator
$(\id_Y\mmb\al)\mmb\id_X$ is {\isc}. But $Y$ is ESR-flat, and, by
Proposition 4, $\al$ is a morphism of left semi-Ruan modules.
Therefore the operator $\id_Y\mmb\al:Y\mmb Z_1\to Y\mmb Z_2$ is
{\isc}. However, this operator is, of course, a morphism of right
modules; moreover, by Proposition 5, it is a morphism of semi-Ruan
modules. It remains to recall that $X$ is also extremely flat. $\tr$

\medskip
The  property of extreme flatness which we introduced is intimately
connected with the question of the extension of bounded morphisms,
descending from the classical Hahn--Banach Theorem.

\medskip
{\bf Definition 4.} A contractive left module $X$ is   {\it
extremely injective with respect to semi-Ruan modules} or, for
short, {\it ESR-injective}, if, for every isometric morphism
$\al:Y\to Z$ of left semi-Ruan modules and an arbitrary bounded
morphism of left modules $\Phi:Y\to X$, there exists a bounded
morphism of left modules $\Psi:Z\to X$ such that the diagram
$$
\xymatrix{ Y \ar[r]^{\al} \ar [d]^\Phi
& Z \ar[dl]^\Psi\\
 X }
$$
\noindent  is commutative  and $\|\Phi\|=\|\Psi\|$. In other words,
every bounded morphism of left modules from $Y$ into $X$ can be
extended, after the identification of $Y$ with a submodule of $Z$,
to a morphism from $Z$ to $X$ with the same norm.

We  define the  `right'  version of this notion in the obvious
symmetric way,.

Finally, by replacing words  `left module'  by  `bimodule'  and also
`semi-Ruan'  by `Ruan', we obtain the definition of a bimodule, {\it
extremely injective with respect to Ruan bimodules} or, for short,
of an {\it ER-injective bimodule}.

\medskip
{\bf Proposition 9.} (i) {\it  Let $X$ be a contractive left or
right normed module. Then it is ESR-flat if and only if its dual
right or, respectively, left module $X^*$} {\it is ESR-injective.

{\rm (ii)} Let $X$ be a contractive bimodule. Then it is ER-flat if
and only if  its dual bimodule $X^*$ is ER-injective.}

\smallskip
$\tl$ Since the argument is parallel in all three cases, we shall
restrict ourselves to  the case of a given left module.

It is obvious that the assertion that $X^*$ is ESR-injective is
equivalent to the following statement: for every isometric morphism
$\al:Y\to Z$ of right semi-Ruan modules, the operator $\al_*:${\bf
h}$_\bb(Z,X^*)\to${\bf h}$_\bb(Y,X^*):\beta\mt\beta\al$, otherwise,
the relevant restriction operator, is strictly co-isometric. (The
latter property means that our operator maps the closed unit ball in
the domain space onto the closed unit ball in the range space).
According to the law of the adjoint associativity, also called the
exponential law (see, e.g.,~\cite[Chapter III, \S3.8]{cml} or
~\cite[Chapter VI, \S3.2]{he2}), the normed space {\bf
h}$_\bb(Y,X^*)$ coincides with the space $(Y\mmb X)^*$ up to the
isometric isomorphism, taking a morphism $\va:Y\to X^*$ to the
functional $f:Y\mmb X\to\co$, well-defined by $f(y\mmb
x)=[\va(y)](x)$. Similarly, {\bf h}$_\bb(Y,X^*)$ is identified with
$(Z\mmb X)^*$. Moreover, one can easily check that we have a
commutative diagram
$$
\xymatrix@C+20pt{ {\bf h}_\bb(Y,X^*) \ar[r]^{\al_*} \ar[d]
& {\bf h}_\bb(Y,X^*) \ar[d]\\
(Z\mmb X)^* \ar[r]^{\al^\bullet} & {(Y\mmb X)^*} },
$$
\noindent where the vertical arrows depict indicated isometric
isomorphisms of normed spaces, and $\al^\bullet$ is the operator
which is  adjoint to $\al\mmb\id_X:Y\mmb X\to Z\mmb X$.
Consequently, the operators $\al_*$ and $\al^\bullet$ are
simultaneously strictly co-isometric or not. But, as an obvious
corollary (in fact, an equivalent formulation) of the Hahn--Banach
theorem, an adjoint operator is strictly co-isometric if  and only
if the original operator is {\isc}. The rest is clear. $\tr$

\medskip
As a byproduct, we have the following result.

\medskip
{\bf Proposition 10.} {\it Suppose that $X$ is a contractive left,
right or two-sided module, and $X_0$ is a  dense submodule of the
respective type. Then $X$ is ESR- (or, according to the sense, ER-)
flat if and only if the same is true of  $X_0$.}

\smallskip
$\tl$ Indeed, the dual (bi)modules of $X$ and $X_0$ coincide, and
hence they are simultaneously extremely injective or not. Then the
previous proposition works. $\tr$\smallskip

In what follows, an assertion that (bi)modules of this or that class
are ESR-(or ER-) injective, will be refered as a `theorem of the
Arveson-Wittstock type'. This is because assertions of that type
have their origin in the  `genuine'  Arveson-Wittstock theorem of
quantum functional analysis (= operator space theory).
As to   Proposition 9, it suggests a certain way to establish such
theorems, reducing questions  about   extreme injectivity to those
about   extreme flatness.
\bigskip

\centerline{\bf $3$. Extreme flatness of certain modules}

\bigskip

\noindent Choose, in addition to  our canonical Hilbert space $L$,
an arbitrary Hilbert space $H$. In this section we shall prove the
extreme flatness of some (bi)modules  connected with this space.

At first, take the algebraic tensor product $L\ot H$ as a subspace
of $L\mh H$ with the induced norm. It is obviously a submodule with
respect to the outer multiplication in the latter space, considered
in the Example 2.

In what follows, the symbol ${\cal S}(\,\cd\,,\,\cd\,)$ denotes the
space of Schmidt operators between two Hilbert spaces, equipped with
the Schmidt norm $\|a\|_\SS:=tr(a^*a)^{{1}/{2}}\,,$ whereas ${\cal
F}_\SS(\,\cd\,,\,\cd\,)$ denotes its dense normed subspace
consisting of finite-rank operators (that is ${\cal
F}(\,\cd\,,\,\cd\,)$, considered with the Schmidt norm).
 We recall that $L\mh H$, as a normed space, can be identified with the space
${\cal S}(H^c,L)$ by means of the {\isc} {\ism} uniquely defined by
taking the elementary tensor $\xi\ot\eta$ to the rank-one operator
$\xi\bigcirc\eta$ (cf., e.g.,~\cite[Chapter 3, \S4.3]{he1}).
Clearly, this {\isc} {\ism} identifies $L\ot H$ with
$\k_\SS(H^{c},L)$.

Note that the space ${\cal S}(H^c,L)$ is a left contractive module
with respect to the usual operator composition: for $a\in\bb$ and
$b\in{\cal S}(H^c,L)$,
 we set $a\cd b:=ab$. Now, returning to the mentioned {\isc} {\ism}, we see that it actually provides the
 identification of $L\mh H$ with ${\cal S}(H^c,L)$ and of $L\ot H$ with $\k_\SS(H^{c},L)$ as
 left contractive modules. This can be immediately checked on elementary tensors.

 From now on we denote the left module ${\cal F}_\SS(H^c,L)$ briefly by $\X$. Take an arbitrary right
semi-Ruan module $Y$. For a time, the main object of our study will
be the normed space $Y\mmb\X$.

Let $c:L\to H^{c}$ be a bounded operator. Consider the bilinear
operator $${\cal T}^Y_c:Y\times\X\to Y: (y,b)\mt y\cd(bc)\,.
$$ Of course, ${\cal T}^Y_c$ is bounded, and $\|{\cal T}_c^Y\|\le\|c\|$. Furthermore, one can
immediately check that this bilinear operator is balanced. Therefore
(see Section 0), it gives rise to the bounded operator from
$Y\mmb\X$ into $Y$, uniquely defined by $$y\mmb b\mt y\cd(bc);y\in
Y,b\in\X$$ and having   norm $\le\|c\|$. Denote this operator by
$T_c^Y$.

\medskip
{\bf Proposition 11.} {\it Let $u\in Y\mmb\X$ be represented as an
elementary tensor $y\mmb b$} (cf. Proposition 1). {\it Further, let
$P\in\k$ be the projection on ${\rm Im}(b)$. Then $u=y\cd P\mmb b$,
and there exists an operator $c\in\k(L,H^c)$ such that
$T_c^Y(u)=y\cd P$.}

\smallskip
$\tl$ Since, of course, we have $Pb=b$,   formulae (3) give the
first of the desired equalities. Further, it is clear from the fact
that $\dim({\rm Im}(b))<\ii$ that there exists $c\in\k(L,H^c)$ such
that $bc=P$.

The second desired equality  follows immediately. $\tr$

\medskip
{\bf Remark 3.} From this, as a  first application, one can easily
obtain that our normed {\tp} $Y\mmb\X$ coincides with the {\it
algebraic} {\tp} of $Y$ and $\X$ over $\bb$. In other words, the
subspace
$$N_1:={\rm span}\{ x\cd a\ot y-x\ot a\cd y\}$$ is closed in $Y\mmp\X$ (cf. Section 0), and thus the quotient
semi-norm on $(Y\mmp\X)/N_1$ is actually a norm. But we do not need
this observation.

\medskip
Now let $\al:Y\to Z$ be an arbitrary bounded morphism of contractive
right semi-Ruan modules. Then, by virtue of the functorial
properties of the module {\tp} (see Section 0), the operator
$\al\mmb\id_\X:Y\mmb\X\to Z\mmb\X$ appears.

Note that for every $c\in\bb(L,H^c)$ we have the commutative diagram
$$
{\xymatrix@C+20pt{ Y\mmb\X \ar[r]^{T_c^Y} \ar[d]_{\al\mmb\id_\X}
& Y \ar[d]^{{\al}}\\
Z\mmb\X \ar[r]^{T_c^Z} & Z \,.}}\qq\qq (5)
$$
\noindent This can be immediately verified on elementary tensors in
$Y\mmb\X$.

\medskip
{\bf Proposition 12.} {\it If $\al$ is an injective map, then the
same is true of  $\al\mmb\id_\X$.}

\smallskip
$\tl$ Suppose that, for $u\in Y\mmb\X$, we have
$\al\mmb\id_\X(u)=0$. Take $y,P$ and $c$ as in Proposition 11. Then
the commutative diagram above gives $y\cd P=T_c^Y(u)=0$. But this,
of course, means that $u=0$. $\tr$

\medskip
At last, we are ready to prove our main theorem.

\medskip
{\bf Theorem 1.} {\it Let $H$ and $K$ be arbitrary Hilbert spaces.
Then the left contractive modules $L\ot H$ and $L\mh H$ are
ESR-flat.}
%
%
%
%

\smallskip
$\tl$ Taking into account Proposition 10, it is sufficient to show
that the module $\X:= \k_\SS(H^{c},L)$, that is, as we remember,
$L\ot H$ in disguise, have the desired property.

Let $\al:Y\to Z$ be an {\isc} morphism of left modules. Consequently
(cf. Section 0), $\al\mmb\id_\X$ is a contractive operator.
Therefore our task is to prove that, for every $v\in Y\mmb\X$ and
$u:=(\al\mmb\id_\X)(v)$, we have $\|v\|\le\|u\|$.

Take the representation of $u$ as $z\mmb b$, as provided by
Proposition 1 (with $Z$ in the role of $Y$). After this, take the
respective $P$ and $c$, indicated in Proposition 11. Then the
commutative diagram (5) gives $$z\cd
P=T_c^Z(u)=T_c^Z(\al\mmb\id_\X)(v)=\al(y)\,,$$ where $y:=T_c^Y(v)\in
Y$.
 From this we have that $(\al\mmb\id_\X)(y\mmb b)=u$, and, because of Proposition 12,
$v=y\mmb b$. Now, remembering that $\al$ is an {\isc} operator, we
obtain the estimate
$$
\|v\|\le\|y\|\|b\|=\|z\cd P\|\|b\|\le\|z\|\|b\|.
$$
Further, $L\ot H$ is a semi-Ruan module, and hence the same is true
of its  `alter ego'  $\X$. It remains to take the infimum of numbers
$\|z\|\|b\|$ over all possible representations of $u$ as elementary
tensors in the previous estimate, and then to apply Proposition 5.
$\tr$

{\bf Remark 4.} As a matter of fact, every semi-Ruan module is
ESR-flat. This was shown by the referee of our paper in his report.
The argument, which is much more lengthy and sophisticated than the
proof of Theorem 1, is suggested by some results of
Lambert~\cite{lam}.

\medskip
As an immediate corollary of Theorem 1, we have the following
theorem.

\medskip
{\bf Theorem 2.} {\it Let $H$ and $K$ be as above. Then:\smallskip

 {\rm (i)} the right contractive modules $L^{c}\ot H$ and $L^{c}\mh H$ are ESR-flat;\smallskip

 {\rm (ii)} the contractive bimodules $(L\ot H)\mmp(L^{c}\ot K)$, $(L\mh H)\mmp(L^{c}\mh K)$, and
their completion $(L\mh H)\mmy(L^{c}\mh K)$ are ER-flat.}

\smallskip
$\tl$ (i)  This follows from the previous theorem and Proposition 7,
being applied to $L\ot H^{c}$ and $L\mh H^{c}$.

(ii)  This follows from the previous theorem, combined with the
assertion (i), Proposition 8 and also, in the case of the third
indicated bimodule, with Proposition 10. $\tr$

\medskip
{\bf Remark 5.} We do not know whether the indicated (bi)modules are
extremely flat in the `absolute' sense. By this, in the case, say,
of left modules, we mean the following property of a given $X$: the
operator $\al\mmb\id_X$ is isometric whenever $\al$ is an {\isc}
morphism between arbitrary (and not only semi-Ruan) right normed
modules.
It is somehow doubtful that the answer is `yes'. Anyhow, if we
consider the similarly defined `absolute extreme' version of
projectivity, then the module $L\mh H$ certainly does not possess
this stronger property provided $\dim H=\ii$. As it was shown
in~\cite{he6}, such a module is not projective even in the usual
sense of Banach homology.

\medskip
Now we came to several Arveson--Wittstock type theorems. Here again
we need the law of the adjoint associativity (= exponential law),
now in a  slightly different version. Namely, suppose that $X$ is a
left and $Y$ is a right contractive module. In such a context,
accordingly to what was said in Section 0, $Y^*$ becomes a left
contractive module, and $\bb(X,Y^*), X\mmp Y$, and $(X\mmp Y)^*$
become contractive bimodules. Then $\bb(X,Y^*)$ coincides with
$(X\mmp Y)^*$ up to the isometric bimodule isomorphism which takes
an operator $\va:X\to Y^*$ to the functional $f:X\mmb Y\to\co$,
well-defined by the formula $f(y\mmb x)=[\va(y)](x)$.

\medskip
{\bf Theorem 3.} {\it Let $H$ and $K$ be arbitrary Hilbert spaces.
Then the left contractive module $L\mh H$ and the right contractive
module $L^{c}\mh H$ are ESR-injective, whereas the contractive
bimodule $\bb(L\mh H,L\mh K)$ {\rm (see Example 1)} is
ER-injective.}

\smallskip
$\tl$ To begin with, it is obvious that, up to an {\isc} {\ism} of
modules of the relevant type, $L\mh H=(L^{c}\mh H^{c})^*$ and
$L^{c}\mh H=(L\mh H^{c})^*$.

Furthermore, by virtue of the Riesz representation theorem, the
normed space $\bb(L\mh H,L\mh K)$, that is $\bb(L\mh H,(L^c\mh
K^c)^c)$, can be identified with the normed space $\bb(L\mh
H,(L^c\mh K^c)^{*})$. Recalling that the latter space is also a
contractive bimodule (of the type $\bb(X,Y^*)$; cf.\ above), we
immediately see that actually we have an identification of
contractive bimodules. Finally, the bimodule $\bb(L\mh H,(L^c\mh
K^c)^{*})$ coincides, by the above mentioned law of the adjoint
associativity, with the  module
 $[(L\mh H)\mmp(L^{c}\mh K^{c})]^*$.

Thus all we have to do in all three cases is to combine Theorems 1
and 2 with Proposition 9. $\tr$

\bigskip

\centerline{\bf $4$. The Arveson--Wittstock theorem}

\bigskip

\noindent In the concluding part of the paper we recall the
Arveson--Wittstock theorem and show that it follows from Theorem 3.
Being, so to say, in the air, it must be well known that this
theorem can be easily deduced from the extension theorems for
morphisms of bimodules. Nevertheless, for the completeness of the
picture, we shall present some details.

In what follows, we use the principal definitions of quantum
functional analysis (= operator space theory) in the frame-work of
the non-coordinate approach. The main ideas of such an approach can
be essentially found in the book of Pisier~\cite{pis} and in the
unpublished notes of Barry Johnson. The detailed definitions, in
somewhat different form, are given in~\cite{he4};
these are the amplification of a linear space and of a linear
operator, a quantum space (= abstract operator space) a concrete
{\qs} (= concrete operator space), and, above all, a {\bc} operator
and its {\bc} norm $\|\cd\|_{cb}$.

\medskip
{\bf Arveson--Wittstock Theorem.} {\it Let $E$ be a quantum subspace
of a {\qs} $G$, and let $H$ be an arbitrary Hilbert space. Then
every {\bc} operator $\va$ from $E$ into the concrete {\qs} $\bb(H)$
can be extended to a {\bc} operator $\psi:G\to\bb(H)$ such that
$\|\psi\|_{cb}=\|\va\|_{cb}$.}

\smallskip
$\tl$ Since $\bb(H)$ is concrete, its amplification $\k\ot\bb(H)$ is
identified with a sub-bimodule of $\bb(L\mh H)$. Let $\Phi$ be a
coextension of the amplification $\va_\ii:=\id_\k\ot\va$ of $\va$ to
a morphism into $\bb(L\mh H)$. Then Theorem 3, being considered for
$Y:=\k\ot E$, $Z:=\k\ot G$, and $K:=H$, provides an extension $\Psi$
of $\Phi$ with the same norm.

Observe that the image of $\Psi$ lies in $\k\ot\bb(H)$. Indeed, $\k\ot Z=\\
{\rm span}\{(\xi\bigcirc\eta)z; \xi,\eta\in L,z\in Z\}\,,$ and, by
the equalities (2), we have $\xi\bigcirc\eta=(\xi\bigcirc
e)p(e\bigcirc\eta)$ for every $e\in L;\|e\|=1$ and $p:=e\bigcirc e$.
Therefore, taking into account the fact that $\Psi$ is a morphism of
$\bb$-bimodules, it is sufficient to show that
$\Psi(pz)=p\cd\Psi(pz)\cd p$ belongs to $\k\ot\bb(H)$. But it is
indeed the case, since it is well known that, for every $\tilde
a\in\bb(L\mh H)$ we have $p\cd\tilde a\cd p=p\mh T$ for
$T\in\bb(H)$, well defined by $e\ot T\xi=(p\mh\id_H)[\tilde
a(e\ot\xi)];\xi\in H$.

Thus $\Psi$ has a well-defined corestriction to $\k\ot\bb(H)$. This
corestriction, being a bimodule morphism, obviously has the form
$\psi_\ii:=\id_\k\ot\psi$ for some operator $\psi:G\to \bb(H)$.
Further, $\|\psi\|_{cb}=\|\psi_\ii\|=\|\Psi\|=\|\Phi\|=
\|\va\|_{cb}$. Finally, $\psi_\ii$ is an extension of $\va_\ii$, and
this obviously implies that $\psi$ is an extension of $\va$. $\tr$
\medskip

{\bf Remark 6.} We should like to emphasize that we have deduced
from Theorem 3 the non-coordinate version of the original
Arveson--Wittstock Theorem, concerning just linear completely
bounded operators. What we did not touch, is the later and more
general form of the Arveson--Wittstock Theorem, dealing with {\bc}
morphisms of bimodules over two arbitrary unital $C^*$-algebras.
Different proofs of such a theorem, formulated in various degrees of
generality, can be found in the papers of Wittstock~\cite[Thm.
3.1]{wi2}, Suen~\cite{sue}, Muhly and Na~\cite[Thm. 3.4]{mun},
Pop~\cite[Thm. 2.5]{pop}. Note that it could be shown that the
respective $\bb$-bimodule version of the Arveson--Wittstock Theorem
and our Theorem 3 are equivalent. This is because, as it was
observed by the referee, there exist isometric functors from the
categories of Ruan and semi-Ruan modules into the categories of
operator $\bb$-bimodules and operator $\bb$-modules, respectively.

\ed